\documentclass[12pt]{article}

\def\qed{\quad \vrule height7.5pt width4.17pt depth0pt} 

\newtheorem{theorem}{Theorem}

\newtheorem{lemma}{Lemma}

\begin{document}

\title{The Critical Parameter for the Heat Equation with a Noise 
Term to Blow Up in Finite Time}
\author{Carl Mueller$^1$
\\Dept. of Mathematics\\University of Rochester
\\Rochester, NY  14627
\\E-mail:  cmlr@troi.cc.rochester.edu}
\date{}
\maketitle

\footnotetext[1]{Supported by an NSA grant. 


AMS 1991 {\em subject classifications}
Primary, 60H15; Secondary, 35R60, 35L05.}

\begin{abstract}
Consider the stochastic partial differential equation
\[
u_t=u_{xx}+u^\gamma\dot{W},
\]
where $x\in\mathbf{I}\equiv [0,J]$, $\dot{W}=\dot{W}(t,x)$ is 2-parameter white noise, 
and we assume that the initial function $u(0,x)$ is nonnegative and not 
identically 0.  We impose Dirichlet boundary conditions on  
$u$ in the interval $\mathbf{I}$.  We say that $u$ blows up in finite time,
with positive probability, if there is a random time $T<\infty$ such that 
\[
P\left(\lim_{t\uparrow T}\sup_x u(t,x)=\infty\right)>0.
\]
It was known that if $\gamma<3/2$, then with probability 1, $u$ does not 
blow up in finite time.  It was also known that there is a positive 
probability of finite time blow-up for $\gamma$ sufficiently large.
    
In this paper, we show that if $\gamma>3/2$, then there is a positive 
probability that $u$ blows up in finite time.
\end{abstract}

\section{Introduction}
\label{intro}
\setcounter{equation}{0}

We consider the heat equation with a nonlinear additive noise term.
\begin{eqnarray}
\label{spde}
u_t&=&u_{xx}+u^\gamma\dot{W},\qquad t>0,x\in\mathbf{I}\equiv [0,J]	\\
u(t,0)&=&u(t,J)=0	\nonumber	\\
u(0,x)&=&u_0(x)  \nonumber
\end{eqnarray}
Here, $\dot W=\dot W(t,x)$ is 2-parameter white noise, $\gamma\ge 1$,
and $u_0(x)$ is a continuous nonnegative function on $\mathbf{I}$,
vanishing at the endpoints, but not identically zero.  Suppose that we
are working on a probability space $(\Omega,\mathcal{F},P)$, and fix a
point $\omega\in\Omega$.  If there exists a random time 
$T=T(\omega)<\infty$ such that 
\[
\lim_{t\uparrow T}\sup_{x\in\mathbf{I}}u(t,x)=\infty
\]
then we say that $u$ blows up in finite time (for the point $\omega$).

For deterministic partial differential equations, there is a large
literature about blow-up in finite time.  See \cite{fuj66}, \cite{fm85},
\cite{fm86}, \cite{fk92}, and \cite{ln92} for example.  Suppose that we 
are dealing with the equation
\begin{eqnarray*}
\frac{\partial w(t,x)}{\partial t}&=&\Delta w(t,x)+g(w(t,x)) \\
w(0,x)&=&w_0(x).
\end{eqnarray*}
One basic ideas is the following.  Suppose that $g(w)$ increases faster 
than linearly.  If a high peak forms in the solution $w(t,x)$, then the
term $g(w(t,x))$ will win out over the term $\Delta w(t,x)$, the growth
of the peak will be governed by the ordinary differential equation 
\[
w'(t)=g(w(t)).
\]
We can solve this equation explicitly, and its solutions often blow up in
finite time.

On the other hand, for stochastic partial differential equations (SPDE)
there are very few papers about finite-time blow-up.  Apart from the heat 
equation, the author \cite{mue97} has studied the wave equation
\begin{eqnarray*}
\frac{\partial^2 u_t(x)}{\partial t^2}&=&\Delta u_t(x)+g(u_t(x))
               \dot W(t,x),\qquad t>0, x\in\mathbf{R}  \\
\frac{\partial u_0(x)}{\partial t}&=&h_1(x)   \\
u_0(x)&=&h_0(x).
\end{eqnarray*}
where $g(u)$ grows like $u(log u)^\alpha$ for some $\alpha>1$.  For
$g(u)=u^\alpha$ with $\alpha>1$, one would guess that solutions would
solutions would blow up in finite time.  But finite time blow-up is not
known for any value of $\alpha$. 
Similar techniques were used in \cite{mue93s} to obtain a modulus of
continuity for solutions of the wave equation with noise in higher
dimensions, with correlated Gaussian noise instead of white noise.

There is more precise information about the heat equation with noise.
Suppose that $u$ is a solution to (\ref{spde}).  In \cite{mue91l} it was 
shown that if $\gamma<3/2$, then, with
probability 1, $u$ does not blow up in finite time.  Krylov
\cite{kry96} gave another proof of this fact for a more general class of
equations.  The papers \cite{mue97} and \cite{mue98} are also relevant.
We refer the reader to Pardoux \cite{par93} for this and
other questions about parabolic SPDE.  Returning to the question of
blow-up, it was shown in \cite{ms93} that there exists $\gamma_0>1$
such that if $\gamma>\gamma_0$, then with positive probability, $u$
blows up in finite time.  The argument in \cite{ms93} was not sharp
enough to give the best value of $\gamma_0$, and the question of
whether one could take $\gamma_0=3/2$ was left open.  The main theorem
of this paper answers this question in the affirmative.

\begin{theorem}
\label{t1}
Let $u(t,x)$ satisfy (\ref{spde}), and suppose that $\gamma>3/2$.
Then, with positive probability, $u$ blows up in finite time.
\end{theorem}

Of course, Theorem \ref{t1} does not tell us what happens at
$\gamma=3/2$.  Surprisingly, the proof of Theorem \ref{t1} uses many
of the same ideas as in \cite{ms93}, although in a sharper form.

Now we discuss the rigorous meaning of (\ref{spde}), following the
formalism of Walsh \cite{wal86}, chapter 3.  Before giving details,
we set up some notation.  Let $G(t,x,y)$ be the fundamental 
solution of the heat equation on $\mathbf{I}$.  If $G(t,x)$ is written as a 
function of 2 variables, we let $G(t,x)$ be the fundamental solution of 
the heat equation on $\mathbf{R}$.  In other words
\[
G(t,x)=\frac 1{\sqrt{4\pi t}}\exp\left(-\frac{x^2}{4t}\right).
\]
It is well known that
\[
G(t,x,y)\le G(t,x-y).
\]
We regard (\ref{spde}) as shorthand for the following integral equation. 
\begin{equation}
\label{int}
u(t,x)=\int_{\mathbf{I}}G(t,x,y)u_0(y)dy
	+\int_0^t\int_{\mathbf{I}}G(t-s,x,y)g(u(s,y))W(dyds)
\end{equation}
where the final term in (\ref{int}) is a white noise integral in the
sense of \cite{wal86}, Chapter 2.  Because $g(u)$ is locally
Lipschitz, standard arguments show that (\ref{spde}) has a unique
solution $u(t,x)$ valid up to the time $\sigma_L$ at which $|u(t,x)|$
first reaches the level $L$ for some $x\in\mathbf{I}$.  Similar arguments are
given in \cite{wal86}, Theorem 3.2 and Corollary 3.4, and his reasoning 
easily carries over to our case. Letting $L\to\infty$, we find that 
(\ref{spde}) has a unique solution for $t<\sigma$, where 
$\sigma=\lim_{L\to\infty}\sigma_L$.  If $\sigma<\infty$, one has  
\[
\lim_{t\uparrow\sigma}\sup_{x\in\mathbf{I}}|u(t,x)|=\infty.
\]
Our goal is to show that $\sigma=\infty$ with probability 1.  

More generally, we regard 
\begin{eqnarray*}
v_t&=&v_{xx}+g(v)\dot{W},\qquad t>0,x\in\mathbf{I}	\\
v(t,0)&=&v(t,J)=0		\\
v(0,x)&=&v_0(x)  \nonumber
\end{eqnarray*}
as a shorthand for the following integral equation, which may only be valid 
up to some blow-up time.
\[
v(t,x)=\int_{\mathbf{I}}G(t,x,y)v_0(y)dy
	+\int_0^t\int_{\mathbf{I}}G(t-s,x,y)g(v(s,y))W(dyds)
\]
Lastly, we will always work with the $\sigma$-fields 
$\mathcal{F}_t=\mathcal{F}^W_t$ generated by the white noise up to time 
$t$.  That is, $\mathcal{F}_t$ is the $\sigma$-field generated by the 
random variables $\int_0^t\int_{\mathbf{I}}\phi(s,x)W(dxds)$, where $\phi$ varies 
over all continuous functions on $[0,t]\times\mathbf{I}$.

We now summarize the argument in \cite{ms93}, which is based on the analysis 
of the formation of high peaks.  Such peaks will occur with positive 
probability.  We wish to show that, with positive probability, such peaks 
grow until they blow up in finite time.  If a high peak forms, we rescale 
the equation and divide the mass of the peak into a collection of peaks of 
smaller mass, and these peaks evolve almost independently.  In this way we 
compare the evolution of $u$ to a branching process.  Large peaks are 
regarded as particles in this branching process.  Offspring are peaks which
are higher by some factor.  We show that the expected number of offspring is 
greater than one when $\gamma>3/2$, and thus the branching process survives 
with positive probability, corresponding to blowup in finite time.

Finally, we remark that in (\ref{spde}), we could replace $u^\gamma$ with 
a function $g(u)$ satisfying $g(u)>cu^\gamma$ for some $c>0$.  Then Theorem 
\ref{t1} would still hold, provided $\gamma>3/2$.

\section{Proof of Theorem \ref{t1}}
\label{proof}
\setcounter{equation}{0}

We give a proof by contradiction.  Assume that 
\begin{equation}\label{contradiction}
P(\sigma<\infty)=0.
\end{equation}
We recall Lemma 2.4 of \cite{ms93}.

\begin{lemma}
\label{scaling}
Suppose that $u$ solves (\ref{spde}) up to some $\mathcal{F}^W_t$ stopping 
time $\tau$.  Let ${\bar L}>0$.  If we let
\[
\tilde v(t,x)\equiv {\bar L}^{-1}
u\left(t{\bar L}^{4(1-\gamma)},x{\bar L}^{2(1-\gamma)}\right) 
\qquad t\ge 0, \quad x\in\mathbf{I}{\bar L}^{2(\gamma-1)}
\]
then $\tilde v(t,x)$ solves
\begin{eqnarray*}
\frac{\partial \tilde v}{\partial t} 
&=& \frac{\partial^2 \tilde v}{\partial x^2}
	+ b(\tilde v,\tilde \xi) \dot {\tilde W}	\\
\tilde v(t,0)&=&\tilde v\left(t,J{\bar L}^{2(\gamma-1)}\right)=0  \\
\tilde v(0,\cdot)&=&\tilde v_0
	\qquad t\ge 0, \quad 0\le x\le J{\bar L}^{2(\gamma-1)}
\end{eqnarray*}
up to the $\mathcal{F}^{\tilde W}_t$-stopping time
$\tau {\bar L}^{2(1-\gamma)}$, where $\tilde v_0(x) 
\equiv{\bar L}^{-1} u_0\left(x{\bar L}^{2(1-\gamma)}\right)$ for all 
$0\le x\le J{\bar L}^{2(\gamma-1)}$,  for all $t\ge 0$ and 
$0\le x\le J{\bar L}^{2(\gamma-1)}$, and $\tilde W$ is the white noise on 
$\mathcal{B}(\mathbf{R}_+\times [0,J{\bar L}^{2(\gamma-1)}])$ defined by
\[
\tilde W(A) \equiv {\bar L}^{3(\gamma-1)}\int_{\mathbf{R}_+\times\mathbf{I}} 
\chi_A (t{\bar L}^{4(\gamma-1)},x{\bar L}^{2(\gamma-1)}) W(dt,dx)
\]
for all $A$ in $\mathcal{B}(\mathbf{R}_+\times [0,J{\bar L}^{2(\gamma-1)}])$ with 
finite Lebesgue measure.
\end{lemma}

In \cite{ms93} we fixed $\bar L$, which we called $L$, and took $\gamma$
to be very large.  In the current proof we wish to deal with all
$\gamma>3/2$, so we take ${\bar L}$ as our large parameter.  We will find
that the probability of a peak getting up to level ${\bar L}$ is about 
$p=1/{\bar L}$, from the gambler's ruin problem.  Using Lemma 
\ref{scaling}, we will see that after rescaling, a peak of size $L$ gives 
rise to $N=L^{2(\gamma-1)}$ offspring.  Thus, the expected number of
offspring of our initial peak should be 
\[
pN=(1/L)\cdot L^{2(\gamma-1)}=L^{2\gamma-3}.
\]
If $\gamma>3/2$, then $2\gamma-3>0$ and $pN\to\infty$ as $L\to\infty$.
Of course, the above heuristic calculation will suffer from the rough
estimates we make during the course of the proof.  Our hope is that taking 
$pN$ is large enough will compensate for all of our sloppiness.

We will consider solutions $\bar u(t,x)$ to a slightly more general equation 
than (\ref{spde}).
\begin{eqnarray}
\label{spdebar}
\bar u_t&=&\bar u_{xx}+g(\bar u)\dot{W},\qquad t>0,x\in\mathbf{I}\equiv [0,J] \\
\bar u(t,0)&=&\bar u(t,J)=0	\nonumber	\\
\bar u(0,x)&=&u_0(x)  \nonumber
\end{eqnarray}
where $g:[0,\infty)\to[0,\infty)$ is a locally Lipschitz function satisfying 
$g(0)=0$ and $g(u)\ge u^\gamma$ for $u>0$.  The same argument as for 
(\ref{spde}) gives existence and uniqueness of $\bar u(t,x)$ up to the 
blow-up time for $\bar u$.

Let $\varphi(t,x)=\varphi^{(T)}(t,x)$ be a solution of the
backward heat equation 
\begin{equation}
\label{var}
\varphi_t=-\varphi_{xx}\qquad 0\le t\le T, x\in\mathbf{R}
\end{equation}
with ``final condition''
\[
\varphi(T,x)=\frac{1}{\sqrt{4\pi T}}\exp\left(-\frac{x^2}{4T}\right).
\]
Of course, $\varphi(T,x)$ is the heat kernel evaluated at time $T$,
and therefore, for $0\le t\le T$,
\[
\varphi(t,x)=\frac{1}{\sqrt{4\pi (2T-t)}}
    \exp\left(-\frac{x^2}{4(2T-t)}\right).
\]
Next, a short calculation yields the following lemma.

\begin{lemma}
\label{lower}
Let $\varphi(t,x)$ be as in (\ref{var}).  If $0\le t\le T$, then 
\[
\varphi^{(T)}(t,x)\ge \sqrt{2}\varphi^{(T)}(T,x)
\]
\end{lemma}

\noindent
\textbf{Proof of Lemma \ref{lower}.}
Since $0\le t\le T$, we have that $T/(2T-t)\ge 1/2$, and 
$T^{-1}\ge (2T-t)^{-1}$.  Therefore
\begin{eqnarray*}
\frac{\varphi^{(T)}(t,x)}{\varphi^{(T)}(T,x)}
	&=& \sqrt{\frac{T}{2T-t}}
\exp\left(\frac{x^2}{4}\left[T^{-1}-(2T-t)^{-1}\right]\right)	\\
&\ge& \sqrt 2
\end{eqnarray*}
This proves Lemma \ref{lower}.
\qed

For future use, we also compute the $\mathbf{L}^1$ norm of $\varphi(t,x)^a$,
for $a>0$.  We claim that there exists a constant $C=C(a)>0$, not
depending on $T$, such that
\begin{eqnarray}
\label{l1}
\|\varphi(t,x)^a\|_1 &=& \int_\mathbf{I} \varphi(t,x)^a dx	\\
	&=& (4\pi (2T-t))^{-a/2}\int_\mathbf{I} 
	  \exp\left(-\frac{a x^2}{4(2T-t)}\right)dx  \nonumber\\
	&=& C'(a)(2T-t)^{(1-a)/2}		\nonumber\\
	&\le& C(a)T^{(1-a)/2}.		\nonumber
\end{eqnarray}
Now, let
\[
M(t)=\int_\mathbf{I} \varphi(t,x){\bar u}(t,x)dx,\qquad 0\le t\le T.
\]
and note that $M(t)$ is a continuous $\mathcal{F}_t$ martingale for
$0\le t\le T$.  This assertion was proven in Lemma 2.3 of \cite{ms93}.
One can also check it heuristically, by formally differentiating $M(t)$
and applying (\ref{spde}) and (\ref{var}).  Lemma 2.3 of \cite{ms93} also
states that $0\le t\le T$, $M(t)$ has square variation
\begin{equation}
\label{square}
\langle M\rangle_t
=\int_0^{t}\int_\mathbf{I} g({\bar u}(s,x))\varphi(s,x)^2 dxds
\ge\int_0^{t}\int_\mathbf{I} {\bar u}(s,x)^{2\gamma}\varphi(s,x)^2 dxds.
\end{equation}
Of course, $M(t)=M^T(t)$ implicitly depends on $T$.
We now prove the following lower bound on $\langle M\rangle_t$.

\begin{lemma}
\label{angle}
There exists a constant $C_1>0$, not depending on $T$, such that if
$0\le t\le T$, then
\[
\langle M\rangle_t\ge C_1T^{-1/2}\int_0^{t} M(s)^{2\gamma}ds.
\]
\end{lemma}

\noindent
\textbf{Proof of Lemma \ref{angle}}.

Let
\begin{equation}
\label{def}
a=\frac{2\gamma-2}{2\gamma-1}.
\end{equation}
Note that 
\begin{equation}
\label{add}
\frac{2-a}{2\gamma}+a=1
\end{equation}
and
\begin{equation}
\label{and}
\frac{1-a}{2}\cdot (1-2\gamma)=-\frac{1}{2}.
\end{equation}
Furthermore, for $t$ fixed,
\[
\frac{\varphi(t,x)^a}{\|\varphi(t,x)^a\|_1}
\]
is a probability density over $x\in\mathbf{R}$.  Using Jensen's inequality,
(\ref{add}), (\ref{l1}), and (\ref{and}), we find that
\begin{eqnarray}
\label{jensen}
\lefteqn{\int_\mathbf{I} \bar u(s,x)^{2\gamma}\varphi(s,x)^2 dx}	\\
	&=& \|\varphi(s,x)^a\|_1\int_\mathbf{I} {\bar u}(s,x)^{2\gamma}\varphi(s,x)^{2-a}
		\frac{\varphi(s,x)^a}{\|\varphi(s,x)^a\|_1}dx
			\nonumber\\
	&\ge& \|\varphi(s,x)^a\|_1 
		\left(\int_\mathbf{I} {\bar u}(s,x)\varphi(s,x)^{(2-a)/(2\gamma)}
     \frac{\varphi(s,x)^a}{\|\varphi(s,x)^a\|_1}dx\right)^{2\gamma}
			\nonumber\\
	&=& \|\varphi(s,x)^a\|_1^{1-2\gamma} 
		\left(\int_\mathbf{I} {\bar u}(s,x)\varphi(s,x)dx\right)
			\nonumber\\
	&\ge& \left(C(a)T^{(1-a)/2}\right)^{1-2\gamma} M(s)^{2\gamma}
			\nonumber\\
	&=&	C_1T^{-1/2}M(s)^{2\gamma}.
		\nonumber
\end{eqnarray}
where $C_1=C(a)^{1-2\gamma}$, and $a$ was defined in (\ref{def}).
After integrating (\ref{jensen}) over $s\in [0,t]$ and putting this
together with (\ref{square}), we get Lemma \ref{angle}.
\qed

Using Lemma \ref{angle}, it is possible to compare $M(t)$ to a
time-changed Brownian motion.  In the standard way, the new time scale
is given by $\langle M\rangle_t$.  Let 
\[
T(L)=16C_1^{-2}L^8
\]
and consider the following gambler's ruin problem.  Start with $M(0)=2$.  
Let $\tau=\tau(L)$ be the first time $t$ that $M(t)=1$ or $M(t)=L$.  Using 
the optional sampling theorem in the usual way, we deduce that 
$EM(\tau)=2$, and therefore (if $M(0)=2$),
\begin{equation}
\label{gamb}
P(M(\tau)=L)=\frac{1}{L-1}.
\end{equation}
In fact, we wish to show

\begin{lemma}
\label{l.gamb}
If $T=T(L)=16C_1^{-2}L^8$, then 
\[
P\left(M(\tau\wedge T)=L\right)\ge \frac{1}{2(L-1)}.
\]
\end{lemma}

\noindent
\textbf{Proof of Lemma \ref{l.gamb}}.

The definition of $\tau$ implies that for all $t\in [0,\tau]$, we have 
$M(t)\ge 1$.  Therefore, by Lemma \ref{angle}, if $t\in [0,\tau]$ then 
\[
\langle M\rangle_t\ge C_1T^{-1/2} t.
\]
Now $M(t)$ is a continuous supermartingale, so it follows that $M(t)$ is 
greater than or equal to a time-changed Brownian motion with time scale 
$\langle M\rangle_t$.  In other words, for some Brownian motion $B(t)$, we 
have $M(t)\ge 2+B(\langle M\rangle_t)$.  Therefore, since 
\[
\langle M\rangle_T\ge CT^{-1/2}T=CT^{1/2},
\]
we have
\begin{eqnarray*}
P\left(T<\tau\right)
  &=& P\left(T<\tau\le\sigma,1<M(t)<L\mbox{ for $t\in [0,T]$}\right) \\
&\le& P\left(T<\tau,1<2+B(\langle M\rangle_t)<L
	\mbox{ for $0\le t\le T$}\right)\\
&=& P\left(T<\tau,1<2+B(t)<L
	\mbox{ for $0\le t\le\langle M\rangle_T$}\right)\\
&=& P\left(T<\tau,1<2+B(t)<L\mbox{ for $0\le t
            \le C_1T^{1/2}$}\right) \\ 
&\le& P\left(\sup_{t\in [0,C_1T^{1/2}]}B(t)<L-2\right).
\end{eqnarray*}
Using the reflection principle, we continue with
\begin{eqnarray*}
\lefteqn{P\left(T<\tau\right)}  \\ 
&\le& 1-P\left(\sup_{t\in [0,C_1T^{1/2}]}B(t)\ge L-2\right) \\
&=& 1-2P\left(B(C_1T^{1/2})\ge L-2\right) \\
&=& P\left(|B(C_1T^{1/2})|\le L-2\right) \\
&=& \int_{-(L-2)}^{L-2}(2\pi C_1T^{1/2})^{-1/2}
	\exp\left(-\frac{x^2}{2C_1T^{1/2}}\right)dx	\\
&\le& 2(L-2)(2\pi C_1T^{1/2})^{-1/2}	\\
&\le& C_1^{-1/2}LT^{-1/4}.
\end{eqnarray*}
Therefore, if
\[
T=T(L)=16C_1^{-2}L^8.
\]
then
\[
P\left(T<\tau\right)<\frac{1}{2(L-1)}.
\]
Then, by (\ref{gamb}),
\begin{eqnarray*}
P\left(M(\tau\wedge T)=L\right)
	&\ge& P\left(M(\tau\wedge T)=L,T\ge\tau\right)	\\
&=& P\left(M(\tau)=L,T\ge\tau\right) \\
&=& P\left(M(\tau)=L\right)-P\left(M(\tau)=L,T<\tau\right)	\\
&\ge& P\left(M(\tau)=L\right)-P\left(T<\tau\right)	\\
&\ge& \frac{1}{L-1}-\frac{1}{2(L-1)}	\\
&=& \frac{1}{2(L-1)}.
\end{eqnarray*}
This proves Lemma \ref{l.gamb}.
\qed

Now we continue with the proof of Theorem \ref{t1}.   Lemma \ref{l.gamb}
implies that 
\[
p\ge \frac{1}{2(L-1)}\ge (2L)^{-1}
\]
Using Lemma \ref{scaling}, with ${\bar L}=L/2$, we deduce that 
\[
N\ge K^{-1}(L/2)^{2(\gamma-1)}-1\ge K^{-1}(L/4)^{2(1-\gamma}
\]
if $L$ is large enough.  Thus we find that
\begin{equation}
\label{pn}
pN\ge K^{-1}(4L)^{-1}(L/4)^{2(\gamma-1)}
    =K^{-1}4^{1-2\gamma}L^{2(\gamma-3/2)}>1
\end{equation}
if $\gamma>3/2$ and $L$ is large enough.

To finish the proof, we can apply the same argument as in \cite{ms93},
sections 3 and 4.  Since these arguments carry over, word for word, we
will merely summarize the argument here, and refer the reader to
\cite{ms93} for details.  

First, we need to split up the solution $u$.  For this, we quote Lemma
2.5 of \cite{ms93}.  But first, define
\[
b(x,y) \equiv\sqrt{(x+y)^{2\gamma} - y^{2\gamma}}.	
\]
\begin{lemma}
\label{split} For $t\ge 0$, $x\in\mathbf{I}$, $i=1,2,\dots,N$
consider the $N$ recursively defined equations
\begin{eqnarray}
\label{2.6}
{\partial u^i \over \partial t} &=& {\partial^2 u^i\over \partial x^2}
	+ b\left(u^i,\sum_{j=1}^{i-1}u^j\right) \dot W^i \\	
	u^i(t,0)&=&u^i(t,J)=0	\nonumber\\
	u^i(0,\cdot)&=&u_0^i.	\nonumber	
\end{eqnarray}
where $u^0\equiv 0$ by definition.  Here the $\{W^i\}$'s are independent
white noises and the $\{u^i_0\}$ are
some collection of nonnegative initial functions.  
Let us then define the process
\[
\tilde u(t,\cdot)\equiv \left\{\begin{array}{ll}
\sum_{i=1}^N u^i(t, \cdot) 
&\mbox{for $0\le t<\min\{\sigma(u^i):\, i=1,2,\ldots,N\}$}\\
	\infty &\mbox{otherwise}
		\end{array}\right.
\]
for all $t\ge 0$.  Here, $\sigma(u^i)$ denotes the blow-up time $\sigma$
with respect to $u^i$.  For $t\ge 0$, $x\in\mathbf{I}$, we have that 
$\tilde u$ is a solution of
\begin{eqnarray*}
\frac{\partial \tilde u}{\partial t} 
	&=& \frac{\partial^2 \tilde u}{\partial x^2}
	+ \tilde u^\gamma \dot {\tilde W}	\\ 
	{\tilde u}(t,0)&=&{\tilde u}(t,J)=0\\
	{\tilde u}(0,\cdot)&=&\sum_{i=1}^N u_0^i
\end{eqnarray*}
for some white noise $\tilde W$ which is a linear combination of the
$\{W^i\}$.
\end{lemma}

We use Lemma \ref{split} to split up the solution $u$ into the sum of
solutions $u^i$.  Later, we will further split up the $u^i$.  Section 4
of \cite{ms93} explains how to use Lemma \ref{split} to split up $u$ over
and over again, at a sequence of stopping times. Each of these smaller 
solutions will have a larger noise term than in (\ref{spde}), so the 
corresponding total mass martingales $U^i(t)=\int_\mathbf{I} u^i(t,x)dx$ 
will have 
\[
\langle	U^i\rangle_t\ge \int_0^tU^i(s)^{2\gamma}ds.
\]
We need a way to split up $u$, given that a certain integral is
sufficiently large.  The following lemma is an easy modification of 
Proposition 3.2 of \cite{ms93}.

\begin{lemma}
\label{birth}
Let $E^\infty_+(\bar J)$ denote the class of nonnegative $\mathbf{C}^\infty$
functions on $[0,\bar J]$.  There exists a constant $K>0$ such that the 
following holds.  Let $J>4$ be fixed.  Set $\bar J\equiv J2^{2(\gamma-1)}$.
If $N>0$ is and integer, and $f_0\in E^\infty_+(\bar J)$  satisfies
\[
\int_0^{\bar J} \phi\left(t,x2^{2(1-\gamma)};z_0,J\right) 
	f_0(x) dx > KN,
\]
for some $z_0$ in $[1,J-1]$ and some $0\le t \le 1$, then there are
functions $\{f_i:\, i=1,2,\ldots,N\}\subset E^\infty_+(\bar J)$ such that
\[
f_0 = \sum_{i=1}^N f_i
\]
and for each $i=1,2,\ldots,N$,
\begin{equation}
\label{3.17}
\int_0^{\bar J} \phi(0,x;z_i,\bar J) f_i(x) dx \ge 2. 
\end{equation}
for some $z_i$ in $[1,\bar J-1]$
\end{lemma}
In \cite{ms93}, Lemma \ref{birth} was shown for $N=[2^{2\gamma-3}]$, but the 
proof given there also implies the above result.  

Now we continue the main argument.  We can assume without loss of generality 
that 
\[
\int_\mathbf{I} G(T,x,y)u(0,y)dy\ge 2.
\]
If this condition fails, wait until time 1, when it has a positive
probability of holding.  Now wait until time $T$.  By Lemma \ref{gamb}, 
we have that 
\begin{equation}
\label{p}
P\left(\int_\mathbf{I} G(2T,x,y)u(0,y)dy\ge L\right)\ge\frac{1}{2(L-1)}=p
\end{equation}
Let
\begin{equation}
\label{N}
N=K^{-1}L^{2(\gamma-1)}.
\end{equation}
Now perform the scaling as in Lemma \ref{scaling}, with $\bar L=L/2$.  For 
the scaled function $\tilde v$, we see that 
\[
\int_\mathbf{I} G(2T,x,y)u(0,y)dy\ge L^{2(\gamma-1)}=KN.
\]
Then, Lemma \ref{birth} shows that we can decompose 
\[
u(t,x)=\sum_{i=1}^N f_i(x)
\]
such that for some set of points $\{z_i\}_{i-1}^N$, 
\[
\int_\mathbf{I} G(2T,x,y)f_i(y)dy\ge 2.
\]
We use these $f_i$ as initial conditions for new functions $u^i(t,x)$,
which satisfy (\ref{2.6}), and we call these $u^i(t,x)$ offspring of 
$u(t,x)$.

If  
\[
\int_\mathbf{I} G(2T,x,y)u(0,y)dy<L,
\]
then we say that mass has died. 

Repeating the argument, we find that there is mass alive at stage $k$ if 
the branching process of the $u$'s is alive at stage $k$.  But this is a
Galton-Watson process with expected number of offspring at least
\[
pN=K^{-1}L^{2(\gamma-1)}\frac{1}{2(L-1)}\ge 2^{-1}K^{-1}L^{2(\gamma-3/2)}
\] by (\ref{p}) and (\ref{N}).  Therefore, if 
$\gamma>3/2$ and $L$ is large enough, the expected number of offspring is
at least
\[
pN>1
\] 
and there is a positive probability of survival.  But survival
means that there is mass present at each stage.  This, in turn, means
that $u(t,x)$ blows up in finite time.  Therefore, there is a positive
probability of finite time blow-up. 

But this conclusion contradicts our assumption (\ref{contradiction}) that 
$P(\sigma<\infty)=0$.  Thus, Theorem \ref{t1} is proved.
\qed


\begin{thebibliography}{Mue97b}

\bibitem[FK92]{fk92}
S.~Filippas and R.~Kohn.
\newblock Refined asymptotics for the blowup of $u_t -\delta u = u^p$.
\newblock {\em Comm. Pure Appl. Math.}, 45:821--869, 1992.

\bibitem[FM85]{fm85}
A.~Friedman and B.~McLeod.
\newblock Blow-up of positive solutions of semilinear heat equations.
\newblock {\em Indiana U. Math. J.}, 34:425--447, 1985.

\bibitem[FM86]{fm86}
A.~Friedman and B.~McLeod.
\newblock Blow-up of solutions of nonlinear degenerate parabolic equations.
\newblock {\em Archive Rat. Mech. and Anal.}, 96:55--80, 1986.

\bibitem[Fuj66]{fuj66}
H.~Fujita.
\newblock On the blowing up of solutions of the Cauchy problem for 
$u_t = \delta u + u^{1+\alpha}$.
\newblock {\em J. Fac. Sci. Tokyo Sect. 1A Math.}, 13:109--124, 1966.

\bibitem[Kry94]{kry96}
N.V. Krylov.
\newblock On $L_p$-theory of stochastic partial differential equations in 
the whole space.
\newblock {\em SIAM J. Math Anal.}, 27(2):313--340, 1994.

\bibitem[LN92]{ln92}
T.-Y. Lee and W.-M. Ni.
\newblock Global existence, large time behavior and life span of solutions 
of a semilinear parabolic Cauchy problem.
\newblock {\em Trans. Amer. Math. Soc.}, 33(1):365--378, 1992.

\bibitem[MS93]{ms93}
C.~Mueller and R.~Sowers.
\newblock Blow-up for the heat equation with a noise term.
\newblock {\em Prob. Th. Rel. Fields}, 97:287--320, 1993.

\bibitem[Mue91]{mue91l}
C.~Mueller.
\newblock Long time existence for the heat equation with a noise term.
\newblock {\em Prob. Th. Rel. Fields}, 90:505--518, 1991.

\bibitem[Mue93]{mue93s}
C.~Mueller.
\newblock A modulus for the 3\_dimensional wave equation with noise: dealing
  with a singular kernel.
\newblock {\em Can. J. Math.}, 45(6):1263--1275, 1993.

\bibitem[Mue97a]{mue97}
C.~Mueller.
\newblock Long-time existence for signed solutions to the heat equation with a
  noise term.
\newblock {\em Ann. Prob.}, 24(1):377--398, 1997.

\bibitem[Mue97b]{mue98}
C.~Mueller.
\newblock Long time existence for the wave equation with a noise term.
\newblock {\em Ann. Prob.}, 25(1):133--152, 1997.

\bibitem[Par93]{par93}
E.~Pardoux.
\newblock Stochastic partial differential equations, a review.
\newblock {\em Bull. Sc. Math.}, 117:29--47, 1993.

\bibitem[Wal86]{wal86}
J.B. Walsh.
\newblock An introduction to stochastic partial differential equations.
\newblock In P.~L. Hennequin, editor, {\em Ecole d'Ete de Probabilites de Saint
  Flour XIV-1984, Lecture Notes in Math. 1180}, Berlin, Heidelberg, New York,
  1986. Springer-Verlag.

\end{thebibliography}

\end{document}